\documentclass{amsart}

\usepackage{amssymb}
\usepackage{amsmath}
\usepackage{amsthm}
\usepackage{amsbsy}

\theoremstyle{plain}
\newtheorem{theorem}{Theorem}[section]
\newtheorem{lemma}[theorem]{Lemma}

\theoremstyle{definition}

\theoremstyle{remark}
\newtheorem{remark}[theorem]{Remark}
\newtheorem*{acknowledgements}{Acknowledgements}

\newcommand{\Z}{\mathbb{Z}}
\newcommand{\R}{\mathbb{R}}

\newcommand{\M}{\tilde M}

\newcommand{\F}{\mathbb{F}}

\newcommand{\del}{\partial}

\begin{document} 

\title{Embedded spheres in $S^2\times S^1\#\dots\#S^2\times S^1$}

\author{Siddhartha Gadgil}

\address{	Stat Math Unit,\\
		Indian Statistical Institute,\\
		Bangalore 560059, India}

\email{gadgil@isibang.ac.in}

\date{\today}


\begin{abstract}
We give an algorithm to decide which elements of $\pi_2(\#_k S^2\times
S^1)$ can be represented by embedded spheres. Such spheres correspond
to splittings of the free group $F_k$ on $k$
generators. Equivalently our algorithm decides whether, for a
handlebody $N$, an element in $\pi_2(N,\del N)$ can be represented by
an embedded disc. We also give an algorithm to decide when classes in
$\pi_2(\#_k S^2\times S^1)$ can be represented by disjoint embedded
spheres.

We introduce the splitting complex of a free group which is analogous
to the complex of curves of a surface. We show that the splitting
complex of $\F_k$ embeds in the complex of curves of a surface of
genus $k$ as a quasi-convex subset.

\end{abstract}

\maketitle

\section{Introduction}

We study here embedded spheres in a $3$-manifold of the form $M=\#_k
(S^2\times S^1)$, i.e., the connected sum of $k$ copies of $S^2\times
S^1$. Group theoretically such spheres correspond to splittings of the
free group $\F_k$ on $k$ generators~\cite{St2}. Understanding these is
likely to be useful in studying $Out(\F_k)$, which is the mapping
class group of the $M$, and more generally the mapping class group of
reducible $3$-manifolds.

Splittings of free groups also correspond to properly embedded discs
in handlebodies~\cite{He}\cite{Ja}. Hence all our results can be
rephrased in terms of properly embedded discs in handlebodies. Because
of the relation to Heegaard splittings, our results are likely to be
useful in this formulation. However, for simplicity of notation, we
shall consider spheres in $M$.

The first question we consider is whether a class in $\pi_2(M)$ can be
represented by an embedded sphere in $M$. Let $\tilde M$ be the
universal cover of $M$. Observe that $\pi_2(M)=\pi_2(\M)=H_2(\M)$ by
Hurewicz theorem. We shall implicitly use this identification
throughout. 

We first consider when $A\in H_2(\M)=\pi_2(M)$ can be represented by
an embedded sphere in $\M$. We shall make use of intersection numbers
(and Poincar\'e duality) for non-compact manifolds. Represent $A$ by a
(not necessarily connected) surface in $\M$ (also denoted $A$). Given
a proper map $c:\R\to\M$ which is transversal to $A$, we consider the
algebraic intersection number $c\cdot A$. This depends only on the
homology class of $A$ and the proper homology class of $c$. The
following gives a criterion for $A$ to be represented by an embedded
sphere.

\begin{theorem}\label{univ}
The class $A\in H_2(\M)$ can be represented by an embedded sphere if
and only if for each proper map $c:\R\to \M$, $c\cdot A\in
\{0,1,-1\}$.
\end{theorem}

For an embedded sphere $S\in M$ with lift $\tilde S\in M$, all the
translates of $\tilde S$ are disjoint from $\tilde S$. In particular,
if $A=[\tilde S]$ is the class represented by $\tilde S$, then $A$ and
$gA$ can be represented by disjoint spheres for each deck
transformation $g$. Thus, our next step is to give a criterion for
when two classes $A$ and $B$ in $H_2(\M)$ can be represented by
disjoint spheres.

\begin{theorem}\label{univ2}
Let $A$ and $B$ be classes in $H_2(\M)$ that can be represented by
embedded spheres. Then $A$ and $B$ can be represented by disjoint
embedded spheres if and only if there do not exist proper maps
$c,c':\R\to \M$ with $c\cdot A=1=c\cdot B$ and $c'\cdot A=1=-c'\cdot
B$.
\end{theorem}

The two above theorems let us determine when, for a class $A\in
\pi_2(M)=H_2(\M)$, the homology classes $A$ and $gA$ can be
represented by disjoint spheres for each $g\in \pi_1(M)$. However to
get an embedded sphere in $M$, we need more. Namely, such a sphere $S$
exists if and only if there is a sphere $\tilde S$ disjoint from all
its translates $g\tilde S$.

Our next result shows that this is automatically satisfied.

\begin{theorem}\label{sphr}
Suppose $A\in \pi_2(M)=H_2(\M)$ is a class such that each for each
$g\in\pi_1(M)$, $A$ and $gA$ can be represented by disjoint spheres in
$\M$. Then $A$ can be represented by an embedded sphere $S\in M$.
\end{theorem}

Thus, we have a criterion for deciding which class can be represented
by an embedded sphere. However our criterion \textit{a priori}
involves checking conditions for infinitely many proper maps
$c,c':\R\to \M$ and infinitely many group elements $g$. We shall show that
it suffices to check only finitely many conditions. This gives the
following result.

\begin{theorem}\label{algo}
There is an algorithm that decides whether a class $A\in \pi_2(M)$ can
be represented by an embedded sphere in $M$.
\end{theorem}

Our methods extend to deciding when two classes $A$ and $B$ can be
represented by disjoint spheres in $M$. This is based on an analogue
of Theorem~\ref{sphr}.

\begin{theorem}\label{sphr2}
Suppose $A$ and $B$ are classes in $\pi_2(M)$ that can be represented
by embedded spheres in $M$. Then $A$ and $B$ can be represented by
disjoint spheres in $M$ if and only if for each $g\in \pi_1(M)$, $A$
and $gB$ can be represented by disjoint spheres in $\M$.
\end{theorem}

\begin{theorem}\label{algo2}
There is an algorithm that decides whether classes $A,B\in \pi_2(M)$
can be represented by disjoint embedded spheres in $M$.
\end{theorem}

In group theoretic terms, isotopy classes of embedded spheres in
$M$ correspond to conjugacy classes of splittings of the free
group $\F_k$. Disjoint spheres in $M$ correspond to splittings
compatible up to conjugacy.

We define the \emph{splitting complex} of $\F_k$ in a manner analogous
to the complex of curves, which has proved very useful in the study of
the mapping class group~\cite{Ha}\cite{Iv} as well as $3$-manifold
topology~\cite{Mi}. Namely, we consider a simplicial complex with
vertices corresponding to conjugacy classes of splittings of $\F_k$. A
finite set of vertices bounds a simplex if the corresponding
splittings are compatible up to conjugacy. This gives a simplicial
complex.

We shall see that this is a (quasi-convex) subcomplex of the complex
of curves.

\begin{theorem}
The splitting complex of $\F_k$ is isomorphic to a subcomplex of the
complex of curves of a surface of genus $k$. Further this subcomplex is
a quasi-convex subset of the complex of curves.
\end{theorem}

The construction of the splitting complex can be made for an arbitrary
group. Moreover, we can consider splittings over any class of subgroups,
for example poly-cyclic groups. Indeed the complex of curves is the
splitting complex of a surface group over $\Z$.

\begin{acknowledgements}
We thank G. Ananda Swarup and Dishant Pancholi for helpful conversations.
\end{acknowledgements}

\section{Ends and spheres in $\M$}

We recall the notion of ends of a space. Let $X$ be a topological
space. For a compact set $K\subset X$, let $C(K)$ denote the set of
components of $X-K$. For $L$ compact with $K\subset L$, we have a
natural map $C(L)\to C(K)$. Thus, as compact subsets of $X$ define a
directed system under inclusion, we can define the set of ends $E(X)$
as the inverse limit of the sets $C(K)$.

It is easy to see that a proper map $f:X\to Y$ induces a map $E(X)\to
E(Y)$ and that this is functorial. In particular, the real line $\R$
has two ends which can be regarded as $-\infty$ and $\infty$. Hence a
proper map $c:\R\to X$ gives a pair of ends $c_-$ and $c_+$ of $X$.

Now consider proper maps $c:\R\to \M$. As $\M$ is a union of simply
connected compact sets, the following lemma is straightforward.

\begin{lemma}
There is a one-one correspondence between proper homotopy classes
of maps $c:\R\to \M$ and pairs $(c_-,c_+)\in E(\M)\times E(\M)$
\end{lemma}

We shall refer to a curve $c$ as above as a proper path from $c_-$ to
$c_+$ or as a proper path joining $c_-$ and $c_+$. We denote such a
path $c$ by $(c_-,c_+)$. This is well defined up to proper
homotopy. In particular, for a homology class $A\in H_2(\M)$, the
intersection number $(c_-,c_+)\cdot A$ is well defined and can be
computed using any proper path joining $c_-$ and $c_+$. We shall use
this implicitly throughout.

We now characterise which homology classes in $\M$ can be represented
by embedded spheres.

\begin{proof}[Proof of Theorem~\ref{univ}]

Suppose $A$ can be represented by an embedded sphere $S$. Then the
complement of $S$ consists of two components with closures $X_1$ and
$X_2$. As $S$ is compact, the space of ends of $\M$ is also
partitioned into sets $E_i=E(X_i)$. For a pair of ends $(c_-,c_+)$, if
both $c_-$ and $c_+$ are contained in the same $E_i$, we have a
corresponding proper path $c$ disjoint from $S$. Otherwise we can
choose $c$ intersecting $S$ in one point. In either case, $c\cdot A$
is $0$, $1$ or $-1$. Computing intersection numbers $(c_-,c_+)\cdot A$
using these paths, it follows that $c\cdot A$ is always $0$, $1$ or
$-1$.

Conversely, assume that for each $c=(c_-,c_+)$, $c\cdot A$ is one of
$0$, $1$ or $-1$. Let $A$ be represented by a (not necessarily
connected) smooth, closed surface, which we also denote $A$. Let
$K\supset A$ be a compact, $3$-dimensional, connected manifold
contained in $\M$ such that the closure $W_i$ of each complementary
component of $K$ is non-compact. As $\M$ is simply-connected and $K$
is connected, $N_i=\del W_i$ is connected for each $W_i$. Note that
there are finitely many sets $W_i$ and $E(\M)$ is partitioned into the
sets $E(W_i)$.

We define a relation on the space of ends $E(\M)$ as follows. For a
pair of ends $e_0$ and $e_1$, let $c$ be a proper path joining $e_0$
to $e_1$. We define $e_0\sim e_1$ if $c\cdot A=0$. We shall show that
the relation $\sim$ is an equivalence relation. When $A\neq 0$ we show
that there are exactly two equivalence classes.

We first need a lemma.

\begin{lemma}\label{add}
For ends $e$, $f$ and $g$ of $\M$.
\begin{itemize}
\item $(e,f)\cdot A=-(f,e)\cdot A$
\item $(e,g)\cdot A=(e,f)\cdot A+(f,g)\cdot A$
\end{itemize}
\end{lemma}
\begin{proof}
The first part is immediate from the definitions. Suppose now $e$, $f$
and $g$ are ends and let $c$ and $c'$ be proper paths from $e$ to $f$
and from $f$ to $g$ respectively. Let $k$ be such that $f\in
E(W_k)$. Then there exist $T\in R$ such that $c([T,\infty))\subset
W_k$ and $c'((-\infty,-T])\subset W_k$. Let $\gamma$ be a path in
$W_k$ joining $c(T)$ and $c'(-T)$. Consider the path
$c''=c|_{(-\infty,T]}*\gamma*c'|_{[-T,\infty)}:\R\to \M$. This is a
proper path from $e$ to $g$ and its intersection points with $A$ are
the union of those of $c$ with $A$ and $c'$ with $A$, with the signs
associated to the points of $c''\cap A$ agreeing with the signs for
$c\cap A$ and $c'\cap A$. Computing $(e,g)\cdot A$ using $c''$, we see
$(e,g)\cdot A=(e,f)\cdot A+(f,g)\cdot A$ as claimed.
\end{proof}

By the above, $\sim$ is an equivalence relation. We next show that
there at most most two equivalence classes. This follows from the next
lemma.

\begin{lemma}\label{two}
Suppose $e\not\sim f$ and $e\not\sim g$. Then $f\sim g$ and
$(e,f)\cdot A=(e,g)\cdot A$.
\end{lemma}
\begin{proof}
By Lemma~\ref{add}, we have
$$(f,g)\cdot A= (e,g)\cdot A-(e,f)\cdot A$$ By hypothesis, each of
$(e,g)\cdot A$ and $(e,f)\cdot A$ is $\pm 1$ and their difference
$(f,g)\cdot A$ is $0$, $1$ or $-1$. It follows that $(e,f)\cdot
A=(e,g)\cdot A$ and $(f,g)\cdot A=0$, i.e., $f\sim g$.
\end{proof}

Now as $A\neq 0$ in homology, by Poincar\'e duality there are ends $e$
and $f$ such that $(e,f)\cdot A\neq 0$, i.e., $e\not\sim f$. Thus
there are exactly two equivalence classes of ends which we denote
$E_1$ and $E_2$. 

Next, observe that given two points in $E(W_i)$, for some $i$, there
is a path joining these in the complement of $K$, hence of $A$. It
follows that these are equivalent. Hence each $E(W_i)$ is contained in
$E_1$ or $E_2$. We now construct a proper function
$f:\M\to\R$. Namely, for each $i$, if $E(W_i)\subset E_1$
(respectively $E(W_i)\subset E_2$), we construct a proper function
$f:W_i\to [-1,-\infty)$ (respectively $f:W_i\to [1,\infty)$). We
extend this across $K$ to get a proper function $f:\M\to\R$.

As $\M$ is simply-connected, using standard techniques due to
Whitehead and Stallings~\cite{St1}\cite{St2}, after a proper homotopy
of $f$ we can assume that $S=f^{-1}(0)$ is a sphere. This separates
$\M$ into subsets $X_1$ and $X_2$. By construction,
$E(X_i)=E_i$. Hence by Poincar\'e duality, after possibly changing the
orientation of $S$, $A=[S]$ as claimed.
\end{proof}

\begin{remark}
By construction $S\subset K$.
\end{remark}

\section{Disjoint spheres in $\M$}

Suppose now that $A$ and $B$ are classes in $H_2(\M)=\pi_2(M)$ which
can be represented by embedded spheres $S$ and $T$. We consider the
when $S$ and $T$ can be chosen to be disjoint. Denote the
closures of the components of the complement of $S$ (respectively $T$)
by $X_1$ and $X_2$ (respectively $Y_1$ and $Y_2$) so that $(e,f)\cdot
A=1$ if and only if $e\in X_1$ and $f\in X_2$ and $(e,f)\cdot B=1$ if
and only if $e\in Y_1$ and $f\in Y_2$. Recall that $(f,e)\cdot
A=-(e,f)\cdot A$ and $(f,e)\cdot B=-(e,f)\cdot B$.

Suppose $S$ and $T$ are disjoint. We first consider the case $T\subset
X_2$. Then $X_1$ is contained in one of $Y_1$ and $Y_2$. If
$X_1\subset Y_1$, then for $c'=(e,f)$, if $c'\dot A=1$ then $e\in
X_1\subset Y_1$ hence $(f,e)\cdot B\neq 1$, i.e., $c'\cdot B\neq
-1$. Thus, there does not exist $c'$ with $c'\cdot A=1=-c'\cdot B$.

By considering other cases similarly, we see that there do not exist
proper maps $c,c':\R\to \M$ with $c\cdot A=1=c\cdot B$ and $c'\cdot
A=1=-c'\cdot B$.

Conversely, suppose there do not exist proper maps $c,c':\R\to \M$
with $c\cdot A=1=c\cdot B$ and $c'\cdot A=1=-c'\cdot B$. We define
three equivalence relations $\sim_A$, $\sim_B$ and $\sim$ on
$E(\M)$. Namely, $e\sim_A f$ (respectively $e\sim_B f$) if $(e,f)\cdot
A=0$ (respectively $(e,f)\cdot B=0$) and $e\sim f$ if $e\sim_A f$ and
$e\sim_B f$. We shall see that $\sim$ partitions $E(\M)$ into at most
three equivalence classes.

Let $e\in E(\M)$ be an end. By Lemma~\ref{two}, for ends $f$,
$(e,f)\cdot A$ has only two possible values, $0$ and one of $1$ and
$-1$. By replacing $A$ by $-A$, we assume $(e,f)\cdot A$ is always $0$
or $1$. Similarly, we assume that $(e,f)\cdot B$ is always $0$ or
$1$. Thus, for ends $f$, the pair$ ((e,f)\cdot A, (e,f)\cdot B)$ has
four possible values. By Lemma~\ref{add}, if $((e,f)\cdot A,
(e,f)\cdot B)=((e,g)\cdot A, (e,g)\cdot B)$, then $f\sim g$. Hence
there are at most four equivalences classes under the relation $\sim$.

We need to show that at least one of these classes is empty. If not,
we can find $f$, $g$ and $h$ with $(e,f)\cdot A=1$, $(e,f)\cdot B=0$,
$(e,g)\cdot A=0$, $(e,g)\cdot B=1$, $(e,h)\cdot B=1$ and $(e,h)\cdot
B=1$. Taking $c=(e,h)$ and $c'=(g,f)$, by Lemma~\ref{add} we see that
$c\cdot A=1=c\cdot B$ and $c'\cdot A=1=-c'\cdot B$, a contradiction.

\begin{remark}\label{cpt}
As the four equivalence classes under $\sim$ are the four
intersections $E(X_i)\cap E(Y_j)$, we see that one of these sets must
be empty, i.e. one of the sets $X_i\cap Y_j$ is compact. This is
important in the sequel.
\end{remark}

If $A$ and $B$ are not independent, then either $A=B$ or $A=-B$ as
both $A$ and $B$ are represented by embedded spheres and are hence
primitive. In this case they can be represented by disjoint embedded
spheres. Hence we may assume that they are independent. By Poincar\'e
duality, it follows that there must be three equivalence classes. Let
$e$, $f$ and $g$ represent the equivalence classes. By changing signs
and permuting if necessary, we can assume that $(e,f)\cdot A=1$,
$(e,f)\cdot B=0$, $(e,g)\cdot B=0$ and $(e,g)\cdot B=1$. 

We now proceed as in the previous section. Choose surfaces
representing $A$ and $B$ and a compact submanifold $K$ containing
these as in the previous section. Let $T$ (a tripod) denote the union
of three half lines $R_e$, $R_f$ and $R_g$, each homeomorphic to
$[0,\infty)$, with the points $0$ in all of them identified. We
construct a proper map $f:\M\to T$ by mapping the components of $W_i$
equivalent to $e$ properly onto $T_e$ and analogously for the other
components and extending this over $K$. Let $1_f\in R_f$ and $1_g\in
R_g$ denote points corresponding to $1$. Then using the techniques of
Whitehead and Stallings, after a proper homotopy of $f$,
$S=f^{-1}(1_e)$ and $T=f^{-1}(1_h)$ are disjoint spheres representing
$A$ and $B$.

\section{Intersection numbers and embedded Spheres}

Suppose now that the class $A\in \pi_2(M)=H_2(\M)$ can be represented
by an embedded sphere $S$ in $\M$. Further assume that for all
$g\in\pi_1(M)$, $A$ and $gA$ can be represented by disjoint embedded
spheres. We show that the class $A$ is represented by a splitting of
the free group $G=\pi_1(M)$ and hence an embedded sphere.

This follows from the work of Scott and Swarup~\cite{SS} using
Remark~\ref{cpt}. In our case we only consider splittings over the
trivial group which simplifies our considerations.

Let $X_1$ and $X_2$ be the closures of the complementary components of
$S$. The Cayley graph of $G$ embeds in $\M$ and the vertices can be
identified with elements of $G$. Let $E_i=X_i\cap G$. Then $E_1$ and
$E_2$ form (almost) complementary almost-invariant sets. The
self-intersection number of the set $E_1$ is the number of $g\in G$
such that all the four sets $E_i\cap gE_j$, $i,j\in\{1,2\}$, are
infinite. But by Remark~\ref{cpt}, for each $g\in G$ at least one of
the intersections $X_i\cap gX_j$ is compact which implies that the
corresponding intersection $E_i\cap gE_j$ is finite (as $G$ is a
discrete subset of $\M$). Thus the self-intersection number of $E_1$
is zero.

By a result of Scott and Swarup~\cite{SS}, it follows that there is a
splitting of the group $G$ corresponding to $A$. Hence, by the Knesser
conjecture, there is an embedded sphere representing the class
$A$. This completes the proof of Theorem~\ref{sphr}

The proof of Theorem~\ref{sphr2} is very similar. We use the result of
Scott and Swarup~\cite{SS} that two splittings are compatible if the
intersection number between the corresponding almost invariant sets
vanishes.

\section{The Algorithms}

We now have necessary and sufficient conditions for deciding whether a
class $A\in\pi_2(M)$ can be represented by an embedded sphere in
$M$. However there are \textit{a priori} infinitely many
conditions. To make this into an algorithm, we reduce these to
finitely many conditions.

Firstly, let $\Gamma\subset M$ be a wedge of circles dual to the
spheres in $M=\#_k S^2\times S^1$. Then the universal cover $T$ of
$\Gamma$ is a tree which embeds in $\M$. We observe that the
complement of the spheres lifts to a set in $\tilde M$ with closure
$P$ a fundamental domain. $P$ intersects $T$ in a unique vertex and
each vertex of $T$ is contained in a unique translate of $P$.

Any proper path $c$ is properly homotopic to an edge path in
$T$. Further, $\pi_2(M)=H_2(\tilde M)$ is generated by spheres $S$
which intersect exactly one edge $e$ of $T$ and with $S\cap e$ is a
single point with transversal intersection.

Thus, elements of $\pi_2(M)$ correspond to finite linear combinations
of edges of $T$. Let $A$ be such an element, and let $\tau\subset T$
be a finite subtree containing the support of $A$. Then for an
edge-path $c$, $c\cdot A$ depends only on the finite edge path
$\xi=c\cap \tau$ contained in $\tau$ with endpoints on $\del
\tau$. Further, as $T$ is a tree without any terminal vertices, any
finite edge path $\xi$ in $\tau$ with endpoints on $\del \tau$ is of
the form $\xi=c\cap \tau$ for a proper path $c$. Hence $A$ is
represented by an embedded sphere in $\M$ if and only if for every
finite edge path $\xi$ in $\tau$ with endpoints on $\del \tau$,
$\xi\cdot A$ is $0$, $1$ or $-1$.

Similarly, given two homology classes $A$ and $B$ in $H_2(\M)$, we
have an algorithm to decide whether $A$ and $B$ can be represented by
disjoint embedded spheres by taking $\tau$ containing the supports of both $A$
and $B$.

Finally, if $A$ is a homology class with $\tau$ a tree supporting $A$,
we first verify whether $A$ can be embedded in $\M$. Next there are at
most finitely many elements $g_1$,\dots $g_n$ in $G$ such that
$\tau\cap g_i\tau$ is non-empty. For each of these $g_i$ we check
whether $A$ and $g_iA$ can be represented by disjoint spheres. Assume
henceforth that $A$ has this property.

Let $K$ be the union of the translates of $P$ containing vertices of
$\tau$. Then $K$ is as in the proof of Theorem~\ref{univ}. Thus $A$
can be represented by an embedded sphere $S$ in $K$. If $\tau\cap
g\tau$ is empty, so is $K\cap gK$ and hence $S\cap gS$, i.e. $A$ and
$gA$ can be represented by disjoint embedded spheres. Thus we need to
check only finitely many conditions for finitely many $g_i$, which can
be done algorithmically.

Similar considerations, using Theorem~\ref{sphr2} gives an algorithm
to decide whether two classes in $\pi_2(M)$ (more generally finitely
many classes in $\pi_2(M)$) can be represented by disjoint spheres.

\section{The Splitting Complex}

The \emph{complex of curves} $CS(S)$ of a surface $S$ has proved to be
very useful in studying the mapping class group of a surface as well
as in $3$-manifold topology. We define analogously the \emph{splitting
complex} $SC(\F_k)$ of a free group $\F_k$.

Namely, let $V$ be the set of splittings of $\F_k$ up to conjugacy, or
equivalently the set of properly embedded discs in a handlebody $H_k$ of
genus $k$ up to isotopy. To define a simplicial complex with vertices
$V$, it suffices to specify when a finite subset of $V$, i.e. a finite
collection of splittings, is the set of vertices of a simplex. We
define the splitting complex by specifying that a collection of
splittings bounds a simplex if it is compatible up to conjugacy. In
topological terms vertices corresponding to a collection of embedded
disjoint discs in $H_k$ bound a simplex in $SC(\F_n)$ if they are
isotopic to disjoint embedded discs.

\begin{theorem}
Let $S_k$ be the surface of genus $k$. Then $SC(\F_k)$ is isomorphic
to a connected quasi-convex subcomplex of $CS(S_k)$.
\end{theorem}
\begin{proof}
We interpret $SC(\F_k)$ in terms of discs in $H_k$. Associating to
each disc its boundary gives an embedding of $SC(\F_k)$ in
$CS(\F_k)$. By results of Masur and Minsky~\cite{MM} it follows that
the image is connected and quasi-convex.
\end{proof}

Many fruitful results regarding the complex of curves, in
particular~\cite{Bo}, have been obtained by studying the relation
between distances in the complex of curves and intersection
numbers. Thus one may hope that similar results regarding the
splitting comlex (and hence $Out(\F_k)$) may be obtained using our
methods. A particularly interseting question is hyperbolicity of the
splitting complex.

\bibliographystyle{amsplain}

\begin{thebibliography}{10}

\bibitem{Bo} Bowditch, Brian H.
\textit{Intersection numbers and the hyperbolicity of the complex of curves.}
preprint.

\bibitem{Ha} Harer, John L.  
\textit{Stability of the homology of the
mapping class groups of orientable surfaces.}  
Ann. of Math. (2) \textbf{121} (1985), 215--249.

\bibitem{He} Heil, Wolfgang 
\textit{On Kneser's conjecture for bounded $3$-manifolds.}  
Proc. Cambridge Philos. Soc.  \textbf{71}  (1972), 243--246.

\bibitem{Iv} Ivanov, Nikolai V. 
\textit{Automorphism of complexes of curves and of Teichm\"uller spaces.}  
Internat. Math. Res. Notices \textbf{14} (1997), 651--666.

\bibitem{Ja} Jaco, William 
\textit{Three-manifolds with fundamental group a free product.} 
Bull. Amer. Math. Soc. \textbf{75} (1969) 972--977

\bibitem{MM} Masur, H. A. ; Minsky, Yair.
\textit{Quasiconvexity in the curve complex},
preprint.

\bibitem{Mi} Minsky, Yair.
\textit{The classification of Kleinian surface groups, I: Models and bounds}
preprint.

\bibitem{SS} Scott, Peter; Swarup, Gadde A. 
\textit{Splittings of groups and intersection numbers.}
Geom. Topol.  \textbf{4} (2000), 179--218.

\bibitem{St1} Stallings, John R.
\textit{A topological proof of Grushko's theorem on free products.}
 Math. Z.  \textbf{90}  (1965) 1--8.

\bibitem{St2} Stallings, John
\textit{Group theory and three-dimensional manifolds.}
Yale Mathematical Monographs, \textbf{4},
Yale University Press, New Haven, Conn.-London, 1971.



\end{thebibliography}

\end{document}